\documentclass[11pt,a4paper]{amsart}             
    
\usepackage{graphicx}
\usepackage{amssymb}

\usepackage[english]{babel}
\usepackage{enumerate}

\newtheorem{lemma}{Lemma}[section]

\newtheorem{remark}[lemma]{Remark}
\newtheorem{corollary}[lemma]{Corollary}
\newtheorem{theorem}[lemma]{Theorem}

\title[Coefficient identification $p$-Stokes]{Coefficient identification of the regularized $p$-Stokes equations}
\usepackage{geometry}
\begin{document}
	
	\maketitle
	\begin{center}
	\author{Niko Schmidt\footnote{Kiel University, n\_f\_schmidt@yahoo.de}}
	\end{center}	
	\section*{Abstract}
The Antarctic and Greenland ice sheet simulation is challenging due to unknown parameters in the $p$-Stokes equations. In this work, we prove the existence of a solution to a parameter identification for the ice rheology and the friction coefficient. Additionally, we verify Gâteaux differentiability of the coefficient-to-state operator by extending a similar result for distributed control. Moreover, we have more complicated boundary conditions. We only have to add a small diffusion term and assume the nonlinear exponent, which is given in applications, to be small enough to obtain the results. Finally, we state the adjoint equation and prove existence and uniqueness of a solution for this equation.
	\subsection*{Keywords}
			Parameter identification problem, Gâteaux differentiability, coefficient-to-state operator, glaciology
\section{Introduction}
The simulation of Greenland and Antarctic ice still poses problems, see \cite{Pachauri2014}. Both ice sheets have unknown physical quantities, like the ice rheology or the friction coefficient, see \cite{GilletChaulet2012}, \cite{Larour2012}, and \cite{Morlighem2013}. The known quantity to estimate the ice rheology and the friction coefficient is the surface velocity field of the glacier. 

We formulate the $p$-Stokes equations for the simulation of glaciers with an additional diffusion term, which does not change the solution drastically if we assume some regularity for the solution to the problem without diffusion, see \cite[Theorem $5.9$]{Schmidt2023a}. We formulate a parameter identification problem and verify the existence of a solution to this problem with the additional diffusion. 
The parameter identification problem considers the coefficients of the nonlinear viscosity and the nonlinear sliding, respectively. Next, we verify Gâteaux differentiability of the coefficient-to-state operator. Finally, we formulate the adjoint equation and prove existence and uniqueness of a solution to this equation. The adjoint equation is used to calculate the gradient of the cost function more efficiently. 

As numerical experiments for the parameter identification problem are quite common in the glaciological community, we refer to \cite{Babaniyi2021}, \cite{Petrat2012}, and \cite{Zhao2018}. Note that those formulations differ slightly from our formulation. In particular, glaciological applications do not have the small diffusion term.

The parameter identification problem extends results from Arada's work, see \cite{Arada2012a}. 
Arada considered distributed control with an additional term from the Navier-Stokes equations. Arada only considered Dirichlet boundary conditions. In contrast, we prescribe boundary conditions for a partly sliding glacier. Arada's work was partly extended by considering a space-dependent exponent of the nonlinear term, see \cite{Guerra2013}. Boundary control was, for example, considered for a general class of semilinear equations in \cite{Casas2006}.

Our paper is structured as follows: In chapter $\ref{pStokesIntro}$, we introduce the $p$-Stokes equations. In chapter $\ref{ExistenceOptimalControl}$, we formulate the parameter identification problem and prove existence of an optimal solution. In chapter $\ref{GateauxDiff}$, we verify Gâteaux differentiability of the coefficient-to-state operator. In chapter $\ref{DualEquation}$, we formulate the adjoint equation and prove existence and uniqueness of a solution for this equation. In the final chapter, we give a short conclusion.
\section{The $p$-Stokes equations}\label{pStokesIntro}
Let $N \in \{2,3\}$, $\Omega\subseteq\mathbb{R}^N$ be a Lipschitz domain. Let $\boldsymbol{v}$ be the velocity and
\begin{align*}
(D \boldsymbol{v})_{ij}&=\frac{1}{2}\left(\frac{\partial v_i}{\partial x_j}+\frac{\partial v_j}{\partial x_i}\right)\text{, }i,j\in \{1,...,N\}.
\end{align*}
We define the function $S^p:\mathbb{R}^{N \times N}\to \mathbb{R}^{N \times N}$,
\begin{align*}
S^p(P)=(|P|^2+\delta^2)^{(p-2)/2}P, \quad \delta>0.
\end{align*}
We interpret $S^p$ for vectors $\boldsymbol{v}\in \mathbb{R}^N$ with
\begin{align*}
\boldsymbol{V}:=\left(\begin{array}{ccc}
v_1 & 0 & 0 \\
0 & v_2 & 0 \\
0 & 0 & v_3
\end{array}\right)
\end{align*}
and $S^p(\boldsymbol{v}):=S^p(\boldsymbol{V})$. We write $S^p_{\Omega}$ and $S^p_{\Gamma}$ for the matrix and vector-valued cases because the latter is only used for boundary integrals. Let $\partial \Omega=\Gamma_b\cup \Gamma_d\cup \Gamma_a$ with $|\Gamma_d|>0$. We set for $r>1$
\begin{align*}
V_r:=\{\boldsymbol{v}\in W^{1,r}(\Omega)^N\text{, }\boldsymbol{v}|_{\Gamma_d}=0\text{, }(\boldsymbol{v}\cdot \boldsymbol{n})|_{\Gamma_b}=0\text{, }\mathrm{div}\boldsymbol{v}=0\}
\end{align*}
with norm
\begin{align}
\|\boldsymbol{v}\|_{V_r}^r=\int_{\Omega}\left(\sum_{i,j=1}^N\left(\frac{\partial v_i}{\partial x_j}\right)^2\right)^{r/2}\, d x.
\end{align}
Let $p\in (1,2)$. The strong formulation of the $p$-Stokes equations is finding a velocity $\boldsymbol{v}\in V_p$ and a pressure $\pi \in L^q(\Omega)$ with
\begin{align*}
-\mathrm{div}(BS_{\Omega}^p(D \boldsymbol{v}))+\nabla \pi&=-\rho \boldsymbol{g}\quad \text{on }\Omega
\end{align*}
with the density $\rho$, the gravity acceleration $\boldsymbol{g}$, the ice rheology $B\in L^{\infty}(\Omega)\text{, }B\geq c_1\in (0,\infty)$, and the boundary conditions
\begin{align*}
\boldsymbol{v}&=\boldsymbol{0}\qquad \qquad \quad \text{on }\Gamma_d,
\\
\sigma \cdot \boldsymbol{n}&=\boldsymbol{0}\qquad \qquad \quad \text{on } \Gamma_a,
\\
\sigma_t &= -\tau S_{\Gamma}^s(\boldsymbol{v_t})\quad \text{ on }\Gamma_b,\quad s\in (1,p],
\\
\boldsymbol{v}\cdot \boldsymbol{n}&=0\qquad \qquad \quad \text{on }\Gamma_b,
\end{align*}
with $\tau \in L^{\infty}(\Gamma_b)$, $\tau \geq 0$, the stress tensor $\sigma:=-\pi I+BS_{\Omega}^p(D\boldsymbol{v})$, and the tangential and normal components
\begin{align*}
\boldsymbol{v} = \boldsymbol{v_t}+(\boldsymbol{v}\cdot \boldsymbol{n})\boldsymbol{n},\quad\text{and}\quad \sigma \cdot \boldsymbol{n}=\sigma_t + (\sigma \cdot \boldsymbol{n}\cdot \boldsymbol{n})\boldsymbol{n}.
\end{align*}
To formulate the variational formulation, we define for $\boldsymbol{v},\boldsymbol{\phi}\in V_p$ and the gravity acceleration $\boldsymbol{g}$
\begin{align*}
(\nabla \boldsymbol{v},\nabla \boldsymbol{\phi})&:=\int_{\Omega}\nabla \boldsymbol{v}:\nabla \boldsymbol{\phi}\, d x=\int_{\Omega}\sum_{i,j=1}^N (\nabla \boldsymbol{v})_{ij}\, (\nabla \boldsymbol{\phi})_{ij}\, d x,
\\
(\boldsymbol{v},\boldsymbol{g})&:=\int_{\Omega}\boldsymbol{v}\cdot \boldsymbol{g}\, d x=\int_{\Omega}\sum_{i=1}^N v_ig_i\, d x.
\end{align*}
A function $\boldsymbol{v}\in V_p$ that fulfills
\begin{align}\label{WeakpStokes}
(BS_{\Omega}^p(D\boldsymbol{v}),\nabla \boldsymbol{\phi})+(\tau S_{\Gamma}^s(\boldsymbol{v}),\boldsymbol{\phi})_{\Gamma_b}
=-(\rho \boldsymbol{g},\boldsymbol{\phi})\text{ for all }\boldsymbol{\phi}\in V_p
\end{align}
is a weak solution of the $p$-Stokes equations. Note here that $\boldsymbol{v}\cdot \boldsymbol{n}$ implies $\boldsymbol{v_t}=\boldsymbol{v}$ on $\Gamma_b$. For equation $(\ref{WeakpStokes})$ exists exactly one solution for $\delta=0$, see \cite[Section $4$ Theorem $1$]{Chen2013}. The case $p=2$ is the Stokes problem. In glaciological applications, $p\in (1,2)$ is used. Of special interest is $p\leq 4/3$, see \cite[Section $1.4.3$]{Fowler2021}, and \cite{Bons2018}.
\begin{lemma}\label{IntegrableS} Let $\mu_0,\delta\in [0,\infty)$, $\boldsymbol{v}\in V_p$. Then $S_{\Omega}^p(\boldsymbol{v})\in L^{p'}(\Omega)^{N \times N}$ is fulfilled.
\end{lemma}
\begin{proof} See \cite[Lemma $2.1$]{Schmidt2023a}.
\end{proof}
The function $S_{\Omega}^p$ is strictly monotone and Lipschitz continuous:
\begin{lemma}\label{UpperBound}
	Let $\mu_0,\delta \in [0,\infty)$. We have for all $P,Q\in \mathbb{R}^{N \times N}$
	\begin{align*}
	(S^p_{\Omega}(P)-S^p_{\Omega}(Q)):(P-Q)&\geq c(\delta+|P|+|Q|)^{p-2}|P-Q|^2
	\\
	|S^p_{\Omega}(P)-S^p_{\Omega}(Q)|&\leq C(\delta+|P|+|Q|)^{p-2}|P-Q|
	\end{align*}
\end{lemma}
\begin{proof}
	The proof is in \cite[Lemma 2.4]{Schmidt2023a} for $\mu_0,\delta\in (0,\infty)$ and uses \cite[Lemma $6.3$]{Diening2007}. However, the conditions are also fulfilled for $\mu_0,\delta \in [0,\infty)$.
\end{proof}
\section{Existence of a Tikhonov minimizer}\label{ExistenceOptimalControl}
We introduce an operator $A:V_2\to V_2^*$ with an additional diffusion term $\mu_0>0$:
\begin{align*}
\langle A\boldsymbol{v},\boldsymbol{\phi}\rangle_{V_2^*,V_2}=(BS_{\Omega}(D\boldsymbol{v}),\nabla \boldsymbol{\phi})+(\tau S_{\Gamma}(\boldsymbol{v}),\boldsymbol{\phi})_{\Gamma_b}
+\mu_0(\nabla \boldsymbol{v},\nabla \boldsymbol{\phi})\quad \text{ for all }\boldsymbol{\phi}\in V_2
\end{align*}
with the simplified notation of $S_{\Omega}$ instead of $S_{\Omega}^p$ and $S_{\Gamma}$ instead of $S_{\Gamma}^s$.

We need the diffusion term for the existence of a Tikhonov minimizer and the Lipschitz continuity of the coefficient-to-state operator. A weak solution $\boldsymbol{v}\in V_2$ of the regularized $p$-Stokes equations fulfills
\begin{align*}
\langle A\boldsymbol{v},\boldsymbol{\phi}\rangle_{V_2^*,V_2}=-(\rho \boldsymbol{g},\boldsymbol{\phi})\quad \text{ for all }\boldsymbol{\phi}\in V_2.
\end{align*}
In \cite[Theorem 2.8]{Schmidt2023a}, the existence of exactly one solution for this problem is proved. The functions $B$ and $\tau$ are unknown in glaciological applications, such as the Greenland ice sheet or the Antarctic. We apply lower bounds on $B$ and $\tau$ to verify existence and uniqueness of a solution for the partial differential equation for each $B$ and $\tau$ and upper bounds for existence of a solution for the minimization problem. We need the upper bound for bounding the sequence $(B_k,\tau_k)$ and finding weakly convergent subsequences. Let $c_1,C_1,C_2\in (0,\infty)$ with $c_1< C_1$. We set $W:=\{(B,\tau)\in H^1(\Omega)\times H^1(\Gamma_b);\,c_1\leq B\leq C_1,\, 0\leq \tau\leq C_2\}$. In practical applications, the surface velocity $\boldsymbol{v_t}$ is known on parts $\emptyset\neq \hat{\Gamma}_a\subseteq \Gamma_a$. We assume that $\hat{\Gamma}_a$ is measurable with $|\hat{\Gamma}_a|>0$ and define the projection operator $P:L^2(\Gamma_a)\to L^2(\Gamma_a)$,
\begin{align*}
P(\boldsymbol{v}(x))=\begin{cases}
P(\boldsymbol{v}(x))&, \quad x\in \hat{\Gamma}_a, \\
0 &,\quad x\in \Gamma_a\setminus \hat{\Gamma}_a.
\end{cases}
\end{align*}
A possible projection could be $\boldsymbol{v}\mapsto \boldsymbol{v_t}$ for the tangential component. In practical applications, the regularization terms penalize the gradient of the ice rheology $\nabla B$ or the gradient of the friction coefficient $\nabla \tau$ as it is assumed that both coefficients should not change too fast. We define the minimization problem with cost function $f$ by
\begin{align}\label{OptimalControlProblem}
\begin{split}
\inf_{(B,\tau)\in W}f(B,\tau)&=\frac{1}{2}\|P(tr\boldsymbol{v})-P(\boldsymbol{v_s})\|_{L^2(\Gamma_a)^N}^2+\frac{\epsilon_1}{2}\|\nabla B\|_{L^2(\Omega)^N}^2+\frac{\epsilon_2}{2}\|\nabla \tau\|_{L^2(\Gamma_b)^N}^2,
\\
\text{ s.t. }\langle A_{B,\tau}\boldsymbol{v}(B,\tau),\boldsymbol{\phi}\rangle_{V_2^*,V_2}&=-(\rho \boldsymbol{g},\boldsymbol{\phi})\text{ for all }\boldsymbol{\phi}\in V_2
\end{split}
\end{align}
with the trace operator $tr:H^1(\Omega)^N \to L^2(\partial \Omega)$, $\epsilon_1,\epsilon_2\in (0,\infty)$, and $\boldsymbol{v_s}\in L^2(\Gamma_a)^N$ measurement data.

First, we verify well-posedness of the problem dependent on the exponent $p$. For this purpose, we need Sobolev's embedding theorem:
\begin{lemma}\label{Sobolev}
	Let $\Omega\subseteq \mathbb{R}^N$ be a Lipschitz domain with $N \in \mathbb{N}$. There exists a continuous and compact embedding $H^1(\Omega)=W^{1,2}(\Omega)\hookrightarrow L^q(\Omega)$ for all $1\leq q< \frac{2N}{N-2}$, and $q<\infty$.
\end{lemma}
\begin{proof}
	See \cite[Chapter $8.9$]{Alt2012}.
\end{proof}
Now, we can bound terms in the variational formulation:
\begin{lemma}\label{BoundTerm}
	Let $\mu_0,\delta\in [0,\infty)$, $(B,\tau)\in W$, $\boldsymbol{v},\boldsymbol{\phi}\in V_2$. Then we have
	\begin{align*}
	|(BS_{\Omega}(D\boldsymbol{v}),\nabla \boldsymbol{\phi})|&\leq \|B\|_{L^{2/(2-p)}(\Omega)}\|\boldsymbol{v}\|_{V_2}^{p-1}\|\boldsymbol{\phi}\|_{V_2}<\infty,
	\\
	|(\tau S_{\Gamma}(tr\boldsymbol{v}),\boldsymbol{\phi})|&\leq c\|\tau\|_{L^{2/(2-p)}(\Gamma_b)}\|\boldsymbol{v}\|_{V_2}^{p-1}\|\boldsymbol{\phi}\|_{V_2}<\infty
	\end{align*}
	with $c>0$ for $N=2$. For $N=3$, we additionally assume $p<5/3$.
\end{lemma}
\begin{proof}
	We determine $r$ such that the integral exists. We have for $r\in (2,\infty)$ with the dual exponent $r'\in (1,2)$:
	\begin{align*}
	\left|\int_{\Omega}B S_{\Omega}(D\boldsymbol{v}):\nabla \boldsymbol{\phi}\, d x\right|
	&\leq \|B\|_{L^r(\Omega)}\left(\int_{\Omega}|S_{\Omega}(D\boldsymbol{v})|^{r'}|\nabla \boldsymbol{\phi}|^{r'}\, d x\right)^{1/r'}
	\\
	&\leq
	\|B\|_{L^r(\Omega)}\left(\int_{\Omega}|S_{\Omega}(D\boldsymbol{v})|^{r't'}\, d x\right)^{1/(r't')}\left(\int_{\Omega}|\nabla \boldsymbol{\phi}|^2\, d x\right)^{1/2}
	\end{align*}
	with $t:=2/r'$ and the dual exponent $t'$. Due to
	\begin{align}\label{InequalitySp}
	|S_{\Omega}(D\boldsymbol{v})|=\left|\left(|D\boldsymbol{v}|^2+\delta^2\right)^{(p-2)/2}D\boldsymbol{v}\right|\leq |D\boldsymbol{v}|^{p-1},
	\end{align}
	the second factor is integrable for $r't'(p-1)= 2$. Some calculations using the dual exponent yield $r=2/(2-p)$.

	%
	%
	%
	
	Due to $B\in W^{1,2}(\Omega)$, we formulate conditions for $W^{1,2}(\Omega)\hookrightarrow L^r(\Omega)=L^{2/(2-p)}(\Omega)$ by using Sobolev's embedding theorem, see Lemma $\ref{Sobolev}$: For $N=2$, the conditions are fulfilled for all $p\in (1,2)$. 
	
	For $N=3$, we have
	\begin{align*}
	\frac{2}{2-p}<6 \Leftrightarrow 2<12-6p \Leftrightarrow 6p<10 \Leftrightarrow p<5/3.
	\end{align*}
	The proof for the second claim is nearly identical. We only use $\Gamma_b$ instead of $\Omega$, $s$ instead of $p$, and the trace operator to bound $\boldsymbol{v}$ and $\boldsymbol{\phi}$ in $V_2$ instead of $L^2(\Gamma_b)^N$.
\end{proof}
\begin{corollary}\label{Wellposed} Let $\mu_0\in (0,\infty)$, $\delta \in [0,\infty)$. The constraint in
	problem $(\ref{OptimalControlProblem})$ can be solved uniquely for each $(B,\tau)\in W$.
\end{corollary}
\begin{proof} Lemma $\ref{BoundTerm}$ yields well-posedness of the variational formulation and continuity of the operator $A$. The proof of the coercivity is identical to \cite[Lemma $2.6$]{Schmidt2023a}, and the proof of the strict monotonicity is identical to \cite[Lemma $2.7$]{Schmidt2023a}, where $\mu_0,\delta \in (0,\infty)$ have been proved. Thus, the Browder-Minty Theorem, \cite{Browder1963}, yields a unique solution for each coefficient $(B,\tau)\in W$.
\end{proof}
We prove the existence of a solution to problem $(\ref{OptimalControlProblem})$ in two steps:
\begin{itemize}
	\item Verifying that $\boldsymbol{v}=\boldsymbol{v}(B,\tau)$ can be bounded independently of $(B,\tau)\in W$.
	\item Showing that there exists a Tikhonov minimizer.
\end{itemize}
\begin{lemma}\label{Boundv}
	For $\mu_0\in (0,\infty)$, $\delta\in [0,\infty)$ the unique solution $\boldsymbol{v}=\boldsymbol{v}(B,\tau)$ is bounded in $V_2$ independently of $(B,\tau)\in W$.
\end{lemma}
\begin{proof} We conclude with Lemma $\ref{UpperBound}$ the first of the following inequalities:
	\begin{align*}
	\mu_0 \|\boldsymbol{v}\|_{V_2}^2\leq \langle A\boldsymbol{v},\boldsymbol{v}\rangle_{V_2^*,V_2}=-(\rho \boldsymbol{g},\boldsymbol{v})\leq \|\rho \boldsymbol{g}\|_{L^2(\Omega)^N}\|\boldsymbol{v}\|_{V_2}
	\Leftrightarrow \|\boldsymbol{v}\|_{V_2}\leq \|\rho \boldsymbol{g}\|_{L^2(\Omega)^N}/\mu_0. 
	\end{align*}
\end{proof}
\begin{theorem}\label{ExOptControl} There exists a Tikhonov minimizer for problem $(\ref{OptimalControlProblem})$ for $\mu_0\in (0,\infty)$, $\delta\in [0,\infty)$ for the case $N=2$ with $p\in (1,2)$. For $N=3$, we additionally have to assume $p<5/3$.
\end{theorem}
\begin{proof} (We follow the usual direct method of calculus of variations. The ideas of bounding $(S^p(D\boldsymbol{v_k}))_k$ for a sequence $\boldsymbol{v_k}\in V_2$ and using the strict monotonicity of the $p$-Stokes equations are adapted from \cite[Theorem $4.1$]{Arada2012}.) We verified in Corollary $\ref{Wellposed}$ that the minimization problem $(\ref{OptimalControlProblem})$ is well-defined. Because the problem has a lower bound, the infimum
	\begin{align*}
	f^*:=\inf_{(B,\tau)\in W}f(B,\tau)
	\end{align*}
	exists. Hence, there exists a sequence $(B_k,\tau_k)_k\in W^{\mathbb{N}}$ with $f(B_k,\tau_k)\to f^*$. All following subsequences will be indexed by $k$ to simplify the notations. 
	
	Let $(B_0,\tau_0)\in W$. A solution $(B,\tau)\in W$ has to fulfill $\epsilon_1\|\nabla B\|_{L^2(\Omega)^N}^2/2\leq f(B_0,\tau_0)$ and $B\leq C_2$ almost everywhere. Thus, the sequence $(B_k)_k$ is bounded in $H^1(\Omega)$. Therefore, a subsequence $B_k\rightharpoonup B$ in $H^1(\Omega)$ exists. 
	Lemma $\ref{BoundTerm}$ yields $B_k\to B$ in $L^{2/(2-p)}(\Omega)$. The same arguments yield $\tau_k \to \tau$ in $L^{2/(2-p)}(\Gamma_b)$. 
	The Riesz-Fischer theorem implies that there exist subsequences $(B_k)_k$ and $(\tau_k)_k$ with $B_k\to B$ almost everywhere and $\tau_k\to \tau$ almost everywhere. The sequence $(\boldsymbol{v_k})_k$ is bounded in $V_2$ due to Lemma $\ref{Boundv}$. Hence, there exists a subsequence $(\boldsymbol{v_k})_k$ with $\boldsymbol{v_k}\rightharpoonup \boldsymbol{v}$ in $V_2$. As the trace operator is a compact operator, see, e.g., \cite[A6.13]{Alt2012}, it follows already the strong convergence $tr\boldsymbol{v_k}\to tr\boldsymbol{v}$ in $L^2(\partial \Omega)$.
	
	We prove that $\boldsymbol{v}$ is a solution to the coefficients $(B,\tau)$. With $\boldsymbol{v_k}\rightharpoonup \boldsymbol{v}$ in $V_2$ the boundedness of $(D\boldsymbol{v_k})_k$ in $L^2(\Omega)^{N \times N}$ follows. Lemma $\ref{IntegrableS}$ implies that $(S_{\Omega}(D\boldsymbol{v_k}))_k$ is bounded in $L^{p'}(\Omega)^{N \times N}$.
	Thus, we find a subsequence with $S_{\Omega}(D\boldsymbol{v_k})\rightharpoonup c_{\Omega}$ in $L^{p'}(\Omega)^{N \times N}$. 
	With $tr\boldsymbol{v_k}\to tr\boldsymbol{v}$ in $L^2(\Gamma_b)$, the sequence $(S_{\Gamma}(tr\boldsymbol{v_k}))_k$ is bounded and a subsequence with $S_{\Gamma}(tr\boldsymbol{v_k})\rightharpoonup c_{\Gamma_b}$ in $L^{p'}(\Gamma_b)^N$ exists. 
	We conclude for all $\boldsymbol{\phi}\in V_2$ with the monotonicity of $A_{k}:=A_{{\mu_0},\delta}(B_k,\tau_k)$ and $A_k(\boldsymbol{v_k})=-\rho \boldsymbol{g}\in V_2^*$:
	\begin{align}\label{RelationWithAmu0}
	0&\leq \langle A_k(\boldsymbol{v_k})-A_k(\boldsymbol{\phi}),\boldsymbol{v_k}-\boldsymbol{\phi}\rangle_{V_2^*,V_2} \nonumber
	\\
	&=\langle A_k(\boldsymbol{v_k}),\boldsymbol{v_k}\rangle_{V_2^*,V_2}-\langle A_k(\boldsymbol{v_k}),\boldsymbol{\phi}\rangle_{V_2^*,V_2}\nonumber
	+\langle A_k(\boldsymbol{\phi}),\boldsymbol{\phi}-\boldsymbol{v_k}\rangle_{V_2^*,V_2} \nonumber
	\\
	&=-(\rho \boldsymbol{g},\boldsymbol{v_k})-\langle A_k(\boldsymbol{v_k}),\boldsymbol{\phi}\rangle_{V_2^*,V_2}+\langle A_k(\boldsymbol{\phi}),\boldsymbol{\phi}-\boldsymbol{v_k}\rangle_{V_2^*,V_2}.
	\end{align}
	For the first summand, we trivially have
	\begin{align*}
	(\rho \boldsymbol{g},\boldsymbol{v_k})\to (\rho \boldsymbol{g},\boldsymbol{v})\quad \text{ for }k\to \infty.
	\end{align*}
	We consider the second summand:
	\begin{align}\label{RightHandSide}
	\langle A_k(\boldsymbol{v_k}),\boldsymbol{\phi}\rangle_{V_2^*,V_2}
	=(B_kS_{\Omega}(D\boldsymbol{v_k}),\nabla \boldsymbol{\phi})
	+(\tau_kS_{\Gamma}(\boldsymbol{v_k}),\boldsymbol{\phi})_{\Gamma_b}+\mu_0(\nabla \boldsymbol{v_k},\nabla \boldsymbol{\phi}).
	\end{align}
	Because of $\boldsymbol{v_k}\rightharpoonup \boldsymbol{v}$ in $V_2$, the third summand on the right-hand side of equation $(\ref{RightHandSide})$ converges. Now, we consider the first summand:
	\begin{align}\label{FirstPartSummandmu0}
	(B_kS_{\Omega}(D\boldsymbol{v_k}),\nabla \boldsymbol{\phi}) 
	=\left((B_k-B)S_{\Omega}(D\boldsymbol{v_k}),\nabla \boldsymbol{\phi}\right)+(BS_{\Omega}(D\boldsymbol{v_k}),\nabla \boldsymbol{\phi}).
	\end{align}
	We obtain for the first summand on the right-hand side of equation $(\ref{FirstPartSummandmu0})$ with Lemma $\ref{BoundTerm}$
	\begin{align}\label{InequalityBound}
	\left|\left((B_k-B)S_{\Omega}(D\boldsymbol{v_k}),\nabla \boldsymbol{\phi}\right)\right|
	\leq& \|B_k-B\|_{L^{2/(2-p)}(\Omega)}\|\boldsymbol{v_k}\|_{V_2}^{p-1}\|\boldsymbol{\phi}\|_{V_2}.
	\end{align}	
	Thus, we obtain for the first summand on the right-hand side of equation $(\ref{FirstPartSummandmu0})$
	\begin{align*}
	\lim_{k \to \infty} \left|\left((B_k-B)S_{\Omega}(D\boldsymbol{v_k}),\nabla \boldsymbol{\phi}\right)\right|=0.
	\end{align*}
	Lemma $\ref{IntegrableS}$ and Lemma $\ref{BoundTerm}$ imply $B\nabla \boldsymbol{\phi}\in L^p(\Omega)\cong L^{p'}(\Omega)^*$. Thus, the second summand on the right-hand side of equation $( \ref{FirstPartSummandmu0})$ converges and we conclude
	\begin{align*}
	\lim_{k \to \infty}(B_kS_{\Omega}(D\boldsymbol{v_k}),\nabla \boldsymbol{\phi})= (Bc_{\Omega},\nabla \boldsymbol{\phi}).
	\end{align*}
	For the boundary, the second summand on the right-hand side of equation $(\ref{RightHandSide})$, we have
	\begin{align}\label{BoundarySummands}
	(\tau_k S_{\Gamma}(tr\boldsymbol{v_k}),\boldsymbol{\phi})
	=\left((\tau_k-\tau)S_{\Gamma}(tr\boldsymbol{v_k}),\boldsymbol{\phi}\right)+(\tau S_{\Gamma}(tr\boldsymbol{v_k}),\boldsymbol{\phi}).
	\end{align}
	We apply the same arguments as before and conclude with Lemma $\ref{BoundTerm}$ and $c\in \mathbb{R}$
	\begin{align*}
	|\left((\tau_k-\tau)S_{\Gamma}(tr\boldsymbol{v_k}),\boldsymbol{\phi}\right)|\leq c \|\tau_k-\tau\|_{L^{2/(2-p)}(\Gamma_b)}\|\boldsymbol{v_k}\|_{V_2}^{p-1}\|\boldsymbol{\phi}\|_{V_2}\to 0.
	\end{align*}
	Because we apply the functional $\tau tr \boldsymbol{\phi}\in L^p(\Gamma_b)^N\cong (L^{p'}(\Gamma_b)^N)^*$ with $S_{\Gamma}(tr \boldsymbol{v_k})$, the second summand of equation $(\ref{BoundarySummands})$ converges. We conclude
	\begin{align*}
	\lim_{k \to \infty}(\tau_k S_{\Gamma}(tr\boldsymbol{v_k}),\boldsymbol{\phi})=(\tau c_{\Gamma_b},\boldsymbol{\phi}).
	\end{align*}
	This yields
	\begin{align}\label{Aklim}
	\begin{split}
	\lim_{k\to \infty}\langle A_{k}(\boldsymbol{v_k}),\boldsymbol{\phi}\rangle_{V_2^*,V_2}
	&=\lim_{k\to \infty}(B_kS_{\Omega}(D\boldsymbol{v_k}),\nabla \boldsymbol{\phi})
	+(\tau_k S_{\Gamma}(\boldsymbol{v_k}),\boldsymbol{\phi})_{\Gamma_b}
	+\mu_0(\nabla \boldsymbol{v_k},\nabla \boldsymbol{\phi})
	\\
	&=(Bc_{\Omega},\nabla \boldsymbol{\phi})
	+(\tau c_{\Gamma_b}, \boldsymbol{\phi})_{\Gamma_b}+\mu_0(\nabla \boldsymbol{v},\nabla \boldsymbol{\phi}).
	\end{split}
	\end{align}
	Now, we consider the last summand on the right-hand side of relation $(\ref{RelationWithAmu0})$:
	\begin{align}\label{NewSummand}
	\begin{split}
	&\quad \langle A_{k}(\boldsymbol{\phi}),\boldsymbol{\phi}-\boldsymbol{v_k}\rangle_{V_2^*,V_2}
	\\
	&=\left(B_k S_{\Omega}(D\boldsymbol{\phi}),\nabla (\boldsymbol{\phi}-\boldsymbol{v_k})\right)
	+\left(\tau_kS_{\Gamma}(\boldsymbol{\phi}), (\boldsymbol{\phi}-\boldsymbol{v_k})\right)_{\Gamma_b}+\mu_0\left(\nabla \boldsymbol{\phi},\nabla (\boldsymbol{\phi}-\boldsymbol{v_k})\right).
	\end{split}
	\end{align}
	We apply the same arguments as before for the first summand on the right-hand side of equation $(\ref{NewSummand})$:
	\begin{align}\label{OtherSummand}
	\left(B_kS_{\Omega}(D\boldsymbol{\phi}),\nabla (\boldsymbol{\phi}-\boldsymbol{v_k})\right)
	=\left((B_k-B)S_{\Omega}(D\boldsymbol{\phi}),\nabla (\boldsymbol{\phi}-\boldsymbol{v_k})\right)+\left(BS_{\Omega}(D\boldsymbol{\phi}),\nabla (\boldsymbol{\phi}-\boldsymbol{v_k})\right).
	\end{align}
	The first summand on the right-hand side of equation $(\ref{OtherSummand})$ vanishes by using the same estimate as in the inequality $(\ref{InequalityBound})$ and using the boundedness of $(\boldsymbol{v_k})_k$ in $V_2$. Thus, we obtain
	\begin{align*}
	\left(B_kS_{\Omega}(D\boldsymbol{\phi}),\nabla (\boldsymbol{\phi}-\boldsymbol{v_k})\right)\rightarrow \left(BS_{\Omega}(D\boldsymbol{\phi}),\nabla (\boldsymbol{\phi}-\boldsymbol{v})\right)\text{ for }k\to \infty.
	\end{align*}
	
	The second summand on the right-hand side of equation $(\ref{NewSummand})$ follows analog, and the third summand is clear. In summary, we can calculate the limit on the right-hand side of relation $(\ref{RelationWithAmu0})$:
	\begin{align}\label{SummaryEquation}
	\begin{split}
	0&\leq \lim_{k \to \infty}-(\rho \boldsymbol{g},\boldsymbol{v_k})-\langle A_{k}(\boldsymbol{v_k}),\boldsymbol{\phi}\rangle_{V_2^*,V_2}+\langle A_{k}(\boldsymbol{\phi}),\boldsymbol{\phi}-\boldsymbol{v_k}\rangle_{V_2^*,V_2}
	\\
	&=-(\rho \boldsymbol{g},\boldsymbol{v})-(Bc_{\Omega},\nabla \boldsymbol{\phi})-(\tau c_{\Gamma_b},\boldsymbol{\phi})_{\Gamma_b}-\mu_0(\nabla \boldsymbol{v},\nabla \boldsymbol{\phi})
	\\
	&\quad+(BS_{\Omega}(D\boldsymbol{\phi}),\nabla (\boldsymbol{\phi}-\boldsymbol{v}))
	+(\tau S_{\Gamma}(\boldsymbol{\phi}),\boldsymbol{\phi}-\boldsymbol{v})
	+\mu_0(\nabla \boldsymbol{v},\nabla (\boldsymbol{\phi}-\boldsymbol{v})).
	\end{split}
	\end{align}
	Because $\boldsymbol{v_k}$ is the solution for the operator $A_k$ we have together with equation $(\ref{Aklim})$ for $\boldsymbol{\phi}=\boldsymbol{v}$
	\begin{align}\label{Solution}
	\begin{split}
	-(\rho \boldsymbol{g},\boldsymbol{v})
	=\lim_{k\to \infty}-(\rho \boldsymbol{g},\boldsymbol{v_k})
	&=\lim_{k \to \infty}\langle A_k(\boldsymbol{v_k}),\boldsymbol{v}\rangle_{V_2^*,V_2}
	\\
	&= (Bc_{\Omega},\nabla \boldsymbol{v})+(\tau c_{\Gamma_b},\boldsymbol{v})_{\Gamma_b}+\mu_0 (\nabla \boldsymbol{v},\nabla \boldsymbol{v}).
	\end{split}
	\end{align}
	Thus, we conclude for the right-hand side of equation $(\ref{SummaryEquation})$
	\begin{align*}
	0&\leq-(\rho \boldsymbol{g},\boldsymbol{v})-(Bc_{\Omega},\nabla \boldsymbol{\phi})-(\tau c_{\Gamma_b},\boldsymbol{\phi})_{\Gamma_b}-\mu_0(\nabla \boldsymbol{v},\nabla \boldsymbol{\phi})
	\\
	&\quad+(BS_{\Omega}(D\boldsymbol{\phi}),\nabla (\boldsymbol{\phi}-\boldsymbol{v}))
	+(\tau S_{\Gamma}(\boldsymbol{\phi}),\boldsymbol{\phi}-\boldsymbol{v})
	+\mu_0(\nabla \boldsymbol{v},\nabla (\boldsymbol{\phi}-\boldsymbol{v}))
	\\
	&= (Bc_{\Omega},\nabla \boldsymbol{v})+(\tau c_{\Gamma_b},\boldsymbol{v})_{\Gamma_b}+\mu_0(\nabla \boldsymbol{v},\nabla \boldsymbol{v})
	-(Bc_{\Omega},\nabla \boldsymbol{\phi})-(\tau c_{\Gamma_b},\boldsymbol{\phi})_{\Gamma_b}-\mu_0(\nabla \boldsymbol{v},\nabla \boldsymbol{\phi})
	\\
	& \quad +(BS_{\Omega}(D\boldsymbol{\phi}),\nabla(\boldsymbol{\phi}-\boldsymbol{v}))+(\tau S_{\Gamma}(\boldsymbol{\phi}),\boldsymbol{\phi}-\boldsymbol{v}))+\mu_0(\nabla \boldsymbol{v},\nabla(\boldsymbol{\phi}-\boldsymbol{v}))
	\\
	&=(Bc_{\Omega},\nabla (\boldsymbol{v}-\boldsymbol{\phi}))+(\tau c_{\Gamma_b},\boldsymbol{v}-\boldsymbol{\phi})_{\Gamma_b}+(BS_{\Omega}(D\boldsymbol{\phi}),\nabla(\boldsymbol{\phi}-\boldsymbol{v}))+(\tau S_{\Gamma}(\boldsymbol{\phi}),\boldsymbol{\phi}-\boldsymbol{v}))
	\\
	&=(B(c_{\Omega}-S_{\Omega}(D\boldsymbol{\phi})),\nabla(\boldsymbol{v}-\boldsymbol{\phi}))+(\tau(c_{\Gamma_b}-S_{\Gamma}(\boldsymbol{\phi})),\boldsymbol{v}-\boldsymbol{\phi})_{\Gamma_b}.
	\end{align*}
	We set $\boldsymbol{\phi}:=\boldsymbol{v}+h\boldsymbol{\psi}$ with $\boldsymbol{\psi} \in V_2$, $\boldsymbol{\psi}|_{\Gamma_b}=0$, and $h\in (0,1)$. This yields
	\begin{align*}
	0\leq h\int_{\Omega}B\left(c_{\Omega}-S_{\Omega}\left(D(\boldsymbol{v}+h\boldsymbol{\psi})\right)\right):\nabla \boldsymbol{\psi}\, d x.
	\end{align*}
	Thus, we conclude
	\begin{align*}
	0\leq \int_{\Omega}B\left(c_{\Omega}-S_{\Omega}\left(D(\boldsymbol{v}+h\boldsymbol{\psi})\right)\right):\nabla \boldsymbol{\psi}\, d x\to \int_{\Omega}B\left(c_{\Omega}-S_{\Omega}(D\boldsymbol{v})\right):\nabla \boldsymbol{\psi}\, d x
	\end{align*}
	for $h\to 0$ because we can apply dominated convergence by the estimate
	\begin{align*}
	|S_{\Omega}(D\boldsymbol{v}+h\psi)|=\left(|D(\boldsymbol{v}+h\boldsymbol{\psi})|^2+\delta^2\right)^{(p-2)/2}|D(\boldsymbol{v}+h\boldsymbol{\psi})|
	&\leq |D(\boldsymbol{v}+h\psi)|^{p-1}
	\\
	&\leq (|D\boldsymbol{v}|+|D\boldsymbol{\psi}|)^{p-1}.
	\end{align*}
	Hence, we have
	\begin{align*}
	0\leq \int_{\Omega}B(c_{\Omega}-S_{\Omega}(D\boldsymbol{v})):\nabla \boldsymbol{\psi}\, d x
	\end{align*}
	for all $\boldsymbol{\psi}\in V_2$ with $\boldsymbol{\psi}|_{\partial \Omega}=0$. With $\psi\in V_2\Rightarrow -\psi \in V_2$ follows
	\begin{align*}
	0=\int_{\Omega}B(c_{\Omega}-S_{\Omega}(D\boldsymbol{v})):\nabla \boldsymbol{\psi}\, d x
	\end{align*}
	for all $\boldsymbol{\psi}\in V_2$ with $\boldsymbol{\psi}|_{\partial \Omega}=0$. Thus, we conclude with the same arguments for the boundary
	\begin{align*}
	\int_{\Omega}Bc_{\Omega}:\nabla \boldsymbol{\psi}\, dx=\int_{\Omega}BS_{\Omega}(D\boldsymbol{v}):\nabla \boldsymbol{\psi}\, dx\quad \text{and }\int_{\Gamma_b}\tau c_{\Gamma_b}\cdot \boldsymbol{\psi}\, ds=\int_{\Gamma_b}\tau S_{\Gamma_b}(\boldsymbol{v})\cdot \psi \, ds
	\end{align*}
	for all $\psi \in V_2$. Thereby, we obtain with relation $(\ref{Aklim})$ and $A_k\boldsymbol{v_k}=\rho \boldsymbol{g} \in V_2^*$
	\begin{align*}
	(\rho \boldsymbol{g},\boldsymbol{\phi})=			\lim_{k\to \infty}\langle A_{k}(\boldsymbol{v_k}),\boldsymbol{\phi}\rangle_{V_2^*,V_2}
	&=\lim_{k\to \infty}(B_kS_{\Omega}(D\boldsymbol{v_k}),\nabla \boldsymbol{\phi})
	+(\tau_k S_{\Gamma}(\boldsymbol{v_k}),\boldsymbol{\phi})_{\Gamma_b}
	+\mu_0(\nabla \boldsymbol{v_k},\nabla \boldsymbol{\phi})
	\\
	&=(Bc_{\Omega},\nabla \boldsymbol{\phi})
	+(\tau c_{\Gamma_b}, \boldsymbol{\phi})_{\Gamma_b}+\mu_0(\nabla \boldsymbol{v},\nabla \boldsymbol{\phi})
	\\
	&=(BS_{\Omega}(D\boldsymbol{v}),\nabla \boldsymbol{\phi})+(\tau c_{\Gamma_b},\boldsymbol{\phi})+\mu_0(\nabla \boldsymbol{v},\nabla \boldsymbol{\phi})
	\end{align*}
	Thus, $\boldsymbol{v}$ is the solution corresponding to the coefficients $(B,\tau)$. Hence, the weak lower semicontinuity of all summands yields
	\begin{align*}
	\inf_{(B,\tau)\in W}f(B,\tau)
	&=\lim_{k\to \infty}f(B_k,\tau_k)
	\\
	&=\lim_{k\to \infty}\frac{1}{2}\|P(tr\boldsymbol{v_k})-P(\boldsymbol{v_s})\|_{L^2(\Gamma_s)^N}^2+\frac{\epsilon_1}{2}\|B_k\|_{H^1(\Omega)}^2+\frac{\epsilon_2}{2}\|\tau_k \|_{H^1(\Gamma_b)}^2
	\\
	&\geq 
	\frac{1}{2}\|P(tr\boldsymbol{v})-P(\boldsymbol{v_s})\|_{L^2(\Gamma_s)^N}^2+\frac{\epsilon_1}{2}\|B\|_{H^1(\Omega)}^2+\frac{\epsilon_2}{2}\|\tau \|_{H^1(\Gamma_b)}^2
	= f(B,\tau).
	\end{align*}
	Therefore, the coefficients $(B,\tau)$ are optimal.
\end{proof}
\begin{remark}
	For $\delta=0$, the proof as in Theorem $\ref{ExOptControl}$ is not possible.
\end{remark}
\begin{proof}
	Following the existence and uniqueness result in \cite{Chen2013}, we would have $\boldsymbol{v}\in V_p$. Thus, inequality $(\ref{InequalityBound})$ would not be sufficient. We would need
	\begin{align*}
	|((B_k-B)S_{\Omega}(D\boldsymbol{v_k}),\nabla \boldsymbol{\phi}|\leq \|B_k-B\|_{L^{\infty}(\Omega)}\|\boldsymbol{v_k}\|_{V_p}^{p-1}\|\boldsymbol{\phi}\|_{V_p}.
	\end{align*}
	However, we do not have convergence of $B_k$ to $B$ in $L^{\infty}(\Omega)$.
\end{proof}
\section{Gâteaux differentiable coefficient-to-state operator}\label{GateauxDiff}
In this section, we prove that the coefficient-to-state operator is Gâteaux differentiable. To obtain this result, we need the Gâteaux differentiability of $A$, which we proved in \cite[Theorem $3.3$]{Schmidt2023a}. There, it was shown that the directional derivative for $\boldsymbol{v},\boldsymbol{w}\in V_2$ is given by
\begin{align}\label{Derivative}
\langle A'(\boldsymbol{v})\boldsymbol{w},\boldsymbol{\phi}\rangle_{V_2^*,V_2}
&=\int_{\Omega}BS_{\Omega}'(D\boldsymbol{v})D\boldsymbol{w}\, :D\boldsymbol{\phi}\, d x+\mu_0 \int_{\Omega}\nabla \boldsymbol{w}:\nabla \boldsymbol{\phi}\, d x +\int_{\Gamma_b}\tau S'_{\Gamma}(\boldsymbol{v})\boldsymbol{w}\cdot \boldsymbol{\phi}\, d s
\end{align}
with the Gâteaux derivatives
\begin{align}\label{SOmegaprime}
\begin{split}
\int_{\Omega}BS_{\Omega}'(D\boldsymbol{v})D\boldsymbol{w}\, :D\boldsymbol{\phi}\, d x
&=
\int_{\Omega}(p-2)B\left(|D\boldsymbol{v}|^2+\delta^2\right)^{(p-4)/2}(D\boldsymbol{v}:D\boldsymbol{w})\, (D\boldsymbol{v}:D \boldsymbol{\phi})\, d x
\\
&\quad+\int_{\Omega}B\left(|D\boldsymbol{v}|^2+\delta^2\right)^{(p-2)/2}D\boldsymbol{w}:D \boldsymbol{\phi}\, d x,
\end{split}    
\\
\int_{\Gamma_b}\tau S'_{\Gamma}(\boldsymbol{v})\boldsymbol{w}\cdot \boldsymbol{\phi}\, d s
&= \int_{\Gamma_b}(s-2)\tau \left(|\boldsymbol{v}|^2+\delta^2\right)^{(s-4)/2}(\boldsymbol{v}\cdot \boldsymbol{w})\, (\boldsymbol{v}\cdot \boldsymbol{\phi})\, d s
\\
&\quad +\int_{\Gamma_b}\tau\left(|\boldsymbol{v}|^2+\delta^2\right)^{(s-2)/2}\boldsymbol{w}\cdot \boldsymbol{\phi}\, d s \nonumber
\end{align}
for all $\boldsymbol{\phi}\in V_2$, $B\in L^{\infty}(\Omega)$ with $c_1\in (0,\infty)$, $B\geq c_1$, $\tau \in L^{\infty}(\Gamma_b)$, and $\tau \geq 0$. 

Now, we introduce variables to formulate the steps for proving the Gâteaux differentiability of the coefficient-to-state operator $\mathcal{S}:W\to V_2$ for $\mu_0,\delta\in (0,\infty)$: 

Let $(\overline{B},\overline{\tau}),(B,\tau)\in W$, $(t_k)_k\in \mathbb{R}^{\mathbb{N}}$ with $t_k\to 0$, $\boldsymbol{\overline{v}}:=\mathcal{S}(\overline{B},\overline{\tau})$, $\boldsymbol{v_k}:=\mathcal{S}(\overline{B}+t_kB,\overline{\tau}+t_k\tau)$, and $\boldsymbol{z_k}:=(\boldsymbol{v_k}-\boldsymbol{\overline{v}})/t_k$.
We prove that the coefficient-to-state operator is Gâteaux differentiable with the following steps:
\begin{itemize}
	\item The operator $\mathcal{S}$ is Lipschitz continuous, see Lemma $\ref{ControlToStateContinuous}$.
	\item The difference $\langle (A(\boldsymbol{v_k})-A(\boldsymbol{\overline{v}}))/t_k,\boldsymbol{\phi}\rangle_{V_2^*,V_2}$ converges for all $\boldsymbol{\phi}\in V_2$, see Lemma $\ref{ConvergenceForGateaux}$.
	\item There exists a subsequence $\boldsymbol{z_k}\rightharpoonup \boldsymbol{\overline{z}}$ in $V_2$, which is a solution of the linearized problem, see Corollary $\ref{SolutionLinearProblem}$.
	\item We conclude $\boldsymbol{z_k}\to \boldsymbol{\overline{z}}$ in $V_2$, see Lemma $\ref{StrongConvergence}$.
\end{itemize}
This idea to prove Gâteaux differentiabilty was used in \cite[Theorem 3.1]{Casas1993} and later in \cite{Slawig2005} and \cite{Arada2012}.
\begin{lemma}\label{ControlToStateContinuous}
	Let $\mu_0\in (0,\infty)$, $\delta \in [0,\infty)$. The coefficient-to-state operator $\mathcal{S}:W\to V_2$, $(B,\tau)\mapsto \boldsymbol{v}(B,\tau)$ is Lipschitz continuous for $N=2$. For $N=3$, we additionally need $p<5/3$.
\end{lemma}
\begin{proof} 
	Let $(B,\tau),(\tilde{B},\tilde{\tau})\in W$. Set $\boldsymbol{v}:=\mathcal{S}(B,\tau)\in V_2$ and $\boldsymbol{\tilde{v}}:=\mathcal{S}(\tilde{B},\tilde{\tau})\in V_2$. Let $\boldsymbol{\phi}\in V_2$.	To prove the Lipschitz continuity, we use Lemma $\ref{UpperBound}$ to obtain the inequality
	\begin{align*}
	&\quad \mu_0\|\boldsymbol{\tilde{v}}-\boldsymbol{v}\|_{V_2}^2
	\\	
	&\leq\mu_0\|\boldsymbol{\tilde{v}}-\boldsymbol{v}\|_{V_2}^2+\left(\tilde{B}\left(S_{\Omega}(D\boldsymbol{\tilde{v}})-S_{\Omega}(D\boldsymbol{v})\right),\nabla (\boldsymbol{\tilde{v}}-\boldsymbol{v})\right)
	+\left(\tilde{\tau} \left(S_{\Gamma}(\boldsymbol{\tilde{v}})-S_{\Gamma}(\boldsymbol{v})\right), (\boldsymbol{\tilde{v}}-\boldsymbol{v})\right)_{\Gamma_b}
	\\
	&=\mu_0\|\boldsymbol{\tilde{v}}-\boldsymbol{v}\|_{V_2}^2	 +\left(\left(\tilde{B}S_{\Omega}(D\boldsymbol{\tilde{v}})-BS_{\Omega}(D\boldsymbol{v})\right),\nabla (\boldsymbol{\tilde{v}}-\boldsymbol{v})\right)
	+\left((B-\tilde{B})S_{\Omega}(D\boldsymbol{v}),\nabla (\boldsymbol{\tilde{v}}-\boldsymbol{v})\right) 
	\\
	&\quad
	+\left((\tilde{\tau} S_{\Gamma}(\boldsymbol{\tilde{v}})-\tau S_{\Gamma}(\boldsymbol{v})), (\boldsymbol{\tilde{v}}-\boldsymbol{v})\right)_{\Gamma_b}
	+\left((\tau-\tilde{\tau})S_{\Gamma}(\boldsymbol{v}), (\boldsymbol{\tilde{v}}-\boldsymbol{v})\right)_{\Gamma_b}
	\\
	&=\langle A_{\tilde{B},\tilde{\tau}}(\boldsymbol{\tilde{v}})-\rho \boldsymbol{g}-A_{B,\tau}(\boldsymbol{v})+\rho \boldsymbol{g},\boldsymbol{\tilde{v}}-\boldsymbol{v}\rangle_{V_2^*,V_2}+\left((B-\tilde{B})S_{\Omega}(D\boldsymbol{v}),\nabla (\boldsymbol{\tilde{v}}-\boldsymbol{v})\right)
	\\
	&\quad+\left((\tau-\tilde{\tau})S_{\Gamma}(\boldsymbol{v}), (\boldsymbol{\tilde{v}}-\boldsymbol{v})\right)_{\Gamma_b}
	\\
	&=\left((B-\tilde{B})S_{\Omega}(D\boldsymbol{v}),\nabla (\boldsymbol{\tilde{v}}-\boldsymbol{v})\right)
	+\left((\tau-\tilde{\tau})S_{\Gamma}(\boldsymbol{v}), (\boldsymbol{\tilde{v}}-\boldsymbol{v})\right).
	\end{align*}
	The last equality followed, because we have for $\boldsymbol{\tilde{v}}$ and $\boldsymbol{v}$
	\begin{align*}
	A_{\tilde{B},\tilde{\tau}}(\boldsymbol{\tilde{v}})-\rho \boldsymbol{g}=0\text{ in }V_2^*\quad \text{and}\quad A_{B,\tau}(\boldsymbol{v})-\rho \boldsymbol{g}=0\text{ in }V_2^*.
	\end{align*}
	Lemma $\ref{BoundTerm}$, and Lemma $\ref{Boundv}$ yield with $\tilde{c}$ independent of $B$, $\tilde{B}$, $\tau$, and $\tilde{\tau}$
	\begin{align*}
	\mu_0 \|\boldsymbol{\tilde{v}}-\boldsymbol{v}\|_{V_2}^2
	&\leq \|B-\tilde{B}\|_{L^r(\Omega)}\|\boldsymbol{v}\|_{V_2}^{p-1}\|\boldsymbol{\tilde{v}}-\boldsymbol{v}\|_{V_2}+c\|\tau-\tilde{\tau}\|_{L^r(\Gamma_b)}\|\boldsymbol{v}\|_{V_2}^{p-1}\|\boldsymbol{\tilde{v}}-\boldsymbol{v}\|_{V_2}
	\\
	&\leq  \tilde{c}\left(\|B-\tilde{B}\|_{L^r(\Omega)}+\|\tau-\tilde{\tau}\|_{L^r(\Gamma_b)}\right)\|\boldsymbol{\tilde{v}}-\boldsymbol{v}\|_{V_2}.
	\end{align*}
	It follows
	\begin{align*}
	\|\boldsymbol{\tilde{v}}-\boldsymbol{v}\|_{V_2}\leq \frac{\tilde{c}}{\mu_0}\left(\|B-\tilde{B}\|_{L^r(\Omega)}+\|\tau-\tilde{\tau}\|_{L^r(\Gamma_b)}\right).
	\end{align*}
	Hence, the coefficient-to-state operator $\mathcal{S}$ is Lipschitz continuous.
\end{proof}
\begin{lemma}\label{ConvergenceForGateaux}
	Let $\mu_0,\delta \in (0,\infty)$. Set $\boldsymbol{z_k}:=(\boldsymbol{v_k}-\boldsymbol{\overline{v}})/t_k$. Let $\boldsymbol{z_k}\rightharpoonup \boldsymbol{\overline{z}}$ in $V_2$ for $t_k\to 0$ and $\boldsymbol{\phi}\in V_2$. Then
	\begin{align*}
	\frac{1}{t_k}\left(\overline{B}(S_{\Omega}(D\boldsymbol{v_k})-S_{\Omega}(D\boldsymbol{\overline{v}})),D \boldsymbol{\phi}\right)&\to 
	(\overline{B}S_{\Omega}'(D\boldsymbol{\overline{v}})D\boldsymbol{\overline{z}},\nabla \boldsymbol{\phi}),
	\\
	\frac{1}{t_k}\left(\overline{\tau}(S_{\Gamma}(\boldsymbol{v_k})-S_{\Gamma}(\boldsymbol{\overline{v}})),\boldsymbol{\phi}\right)&\to 
	(\overline{\tau}S_{\Gamma}'(\boldsymbol{\overline{v}})\boldsymbol{\overline{z}},\boldsymbol{\phi})_{\Gamma_b},
	\\
	\frac{1}{t_k}(\nabla \boldsymbol{v_k}-\nabla \boldsymbol{\overline{v}},\nabla \boldsymbol{\phi})&\to (\nabla \boldsymbol{\overline{z}},\nabla\boldsymbol{\phi}).
	\end{align*}
\end{lemma}
\begin{proof}
	See \cite[Lemma $6.10$]{Arada2012} for the first limit for $3N/(N+2)\leq p\leq 2$. However, the proof for $p\in (1,2)$ is identical. The second limit follows because we only replace $\Omega$ with $\Gamma_b$, use $s$ instead of $p$, and consider $(\boldsymbol{V_k})_{ij}:=I_{ij}(\boldsymbol{v_k})_i$, $\boldsymbol{\overline{V}}_{ij}:=I_{ij}(\boldsymbol{\overline{v}})_i$ with the identity matrix $I$. The third limit follows directly from $\boldsymbol{z_k}\rightharpoonup \boldsymbol{\overline{z}}$ in $V_2$.
\end{proof}
\begin{remark}
	Lemma $\ref{ConvergenceForGateaux}$ is not true for $\delta=0$ as the derivative $S_{\Omega}'$ does not exist for $\delta=0$ in all directions, see \cite[Remark $3.2$]{Schmidt2023a}.
\end{remark}
\begin{corollary}\label{SolutionLinearProblem}
	Let $\mu_0,\delta\in (0,\infty)$ and $(B,\tau),(\overline{B},\overline{\tau})\in W$. Let for $t_k>0$ the functions $\boldsymbol{\overline{v}},\boldsymbol{v_k}\in V_2$ be the solution of
	\begin{align}
	\mu_0(\nabla \boldsymbol{\overline{v}},\nabla \boldsymbol{\phi})+(\overline{B}S_{\Omega}(D\boldsymbol{\overline{v}}),\nabla \boldsymbol{\phi})+(\overline{\tau}S_{\Gamma}(\boldsymbol{v}),\boldsymbol{\phi})_{\Gamma_b})&=-(\rho \boldsymbol{g},\boldsymbol{\phi}),\label{Sol1}
	\\
	\mu_0(\nabla \boldsymbol{v_k},\nabla \boldsymbol{\phi})+\left((\overline{B}+t_kB)S_{\Omega}(D\boldsymbol{v_k}),\nabla \boldsymbol{\phi}\right)+\left((\overline{\tau}+t_k\tau)S_{\Gamma}(\boldsymbol{v_k}),\boldsymbol{\phi}\right)_{\Gamma_b}&=-(\rho \boldsymbol{g},\boldsymbol{\phi})\label{Sol2}
	\end{align}
	for all $\boldsymbol{\phi}\in V_2$.	Set $\boldsymbol{z_k}:=(\boldsymbol{v_k}-\boldsymbol{\overline{v}})/t_k$. Then there exists a subsequence $\boldsymbol{z_k}\rightharpoonup \boldsymbol{\overline{z}}$ in $V_2$ for $t_k \to 0$ with
	\begin{align}\label{Linearity}
	\langle D_{\boldsymbol{v}}A_{\overline{B},\overline{\tau}}(\boldsymbol{\overline{v}})\boldsymbol{\overline{z}},\boldsymbol{\phi}\rangle_{V_2^*,V_2}=-\langle D_{B,\tau}A_{\overline{B},\overline{\tau}}(\boldsymbol{\overline{v}})(B,\tau),\boldsymbol{\phi}\rangle_{V_2^*,V_2}.
	\end{align}
\end{corollary}
\begin{proof}
	We calculate the difference between the equations $(\ref{Sol2})$ and $(\ref{Sol1})$
	\begin{align*}
	&\quad \mu_0(\nabla (\boldsymbol{v_k}-\boldsymbol{\overline{v}}),\nabla \boldsymbol{\phi})+\left((\overline{B}+t_kB)S_{\Omega}(D\boldsymbol{v_k})-\overline{B}S_{\Omega}(D\boldsymbol{\overline{v}}),\nabla \boldsymbol{\phi}\right)
	\\
	&+\left((\overline{\tau}+t_k\tau)S_{\Gamma}(\boldsymbol{v_k})-\overline{\tau}S_{\Gamma}(\boldsymbol{\overline{v}}),\boldsymbol{\phi}\right)_{\Gamma_b}=0.
	\end{align*}
	This is equivalent to
	\begin{align}\label{ConvergenceEquation}
	\begin{split}
	&\quad \mu_0 \left(\frac{\nabla (\boldsymbol{v_k}-\boldsymbol{\overline{v}})}{t_k},\nabla \boldsymbol{\phi}\right)+\left(\overline{B}\frac{S_{\Omega}(D\boldsymbol{v_k})-S_{\Omega}(D\boldsymbol{\overline{v}})}{t_k},\nabla \boldsymbol{\phi}\right)+\left(\overline{\tau}\frac{S_{\Gamma}(\boldsymbol{v_k})-S_{\Gamma}(\boldsymbol{\overline{v}})}{t_k},\boldsymbol{\phi}\right)_{\Gamma_b}
	\\
	&=-(BS_{\Omega}(D\boldsymbol{v_k}),\nabla \boldsymbol{\phi})-(\tau S_{\Gamma}(\boldsymbol{v_k}),\boldsymbol{\phi})_{\Gamma_b}.
	\end{split}
	\end{align}
	Lemma $\ref{ControlToStateContinuous}$ implies that $(\boldsymbol{z_k})_k$ is bounded in $V_2$. Thus, there exists a subsequence with $\boldsymbol{z_k}\rightharpoonup \boldsymbol{\overline{z}}$ in $V_2$.
	Lemma $\ref{ConvergenceForGateaux}$ implies the convergence of the left-hand side of equation $(\ref{ConvergenceEquation})$ to the left-hand side of equation $(\ref{Linearity})$. Due to $\boldsymbol{v_k}\to \boldsymbol{\overline{v}}$ in $V_2$, convergence of the right-hand side follows with dominated convergence. Thus, we obtain the claim.
\end{proof}
Before we can prove strong convergence of $(\boldsymbol{z_k})_k$, we have to verify coercivity of $S'_{\Omega}$ and $S'_{\Gamma}$:
\begin{lemma}\label{LowBound}
	Let $\mu_0,\delta\in (0,\infty)$, $\boldsymbol{v},\boldsymbol{z}\in V_2$, and $(B,\tau)\in W$. Then, we have
	\begin{align*}
	(BS'_{\Omega}(D \boldsymbol{v})D\boldsymbol{z},D\boldsymbol{z})&\geq 0,
	\\
	(\tau S'_{\Gamma}(\boldsymbol{v})\boldsymbol{z},\boldsymbol{z})_{\Gamma_b}&\geq 0.
	\end{align*}
\end{lemma}
\begin{proof}
	We know the first of the following equalities from equation $(\ref{SOmegaprime})$:
	\begin{align*}
	&\quad (B S'_{\Omega}(D\boldsymbol{v})D\boldsymbol{z},D\boldsymbol{z})
	\\
	&=\int_{\Omega}B \left(|D\boldsymbol{v}|^2+\delta^2\right)^{(p-2)/2}D\boldsymbol{z}: D\boldsymbol{z}\, d x 
	\\
	&\quad + \int_{\Omega}B (p-2)\left(|D\boldsymbol{v}|^2+\delta^2\right)^{(p-4)/2}(D\boldsymbol{v}: D\boldsymbol{z})\, (D\boldsymbol{v}: D\boldsymbol{z})\, d x
	\\
	&=\int_{\Omega}B \left(|D\boldsymbol{v}|^2+\delta^2\right)^{(p-2)/2}\left(|D\boldsymbol{z}|^2+(p-2)\frac{(D\boldsymbol{v}\cdot D\boldsymbol{z})\,(D\boldsymbol{v}\cdot D\boldsymbol{z})}{|D\boldsymbol{v}|^2+\delta^2}\right)\, d x.
	\end{align*}
	We use $B>0$ and $p<2$ to conclude with the Cauchy-Schwarz inequality
	\begin{align*}
	&\quad \int_{\Omega}B \left(|D\boldsymbol{v}|^2+\delta^2\right)^{(p-2)/2}\left(|D\boldsymbol{z}|^2+(p-2)\frac{(D\boldsymbol{v}\cdot D\boldsymbol{z})\,(D\boldsymbol{v}\cdot D\boldsymbol{z})}{|D\boldsymbol{v}|^2+\delta^2}\right)\, d x
	\\
	&\geq \int_{\Omega}B \left(|D\boldsymbol{v}|^2+\delta^2\right)^{(p-2)/2}\left(|D\boldsymbol{z}|^2+(p-2)\frac{|D\boldsymbol{v}|^2|D\boldsymbol{z}|^2}{|D\boldsymbol{v}|^2+\delta^2}\right)\, d x
	\\
	&\geq \int_{\Omega}B \left(|D\boldsymbol{v}|^2+\delta^2\right)^{(p-2)/2}\left(|D\boldsymbol{z}|^2+(p-2)|D\boldsymbol{z}|^2\right)\, d x
	\\
	&= \int_{\Omega}B \left(|D\boldsymbol{v}|^2+\delta^2\right)^{(p-2)/2}(p-1)|D\boldsymbol{z}|^2\, d x\geq 0.
	\end{align*}
	The same arguments are valid for $\Gamma_b$ instead of $\Omega$, $s$ instead of $p$, $\boldsymbol{v}$ instead of $D\boldsymbol{v}$, and $\boldsymbol{z}$ instead of $D\boldsymbol{z}$.
\end{proof}
\begin{lemma}\label{StrongConvergence}
	Let the same conditions be fulfilled as in Corollary \ref{SolutionLinearProblem}. Then follows $\boldsymbol{z_k}\to \boldsymbol{\overline{z}}$ in $V_2$.
\end{lemma}
\begin{proof}
	(We transfer the proof in \cite[Proposition $6.11$]{Arada2012} by having more general boundary conditions and having a parameter identification problem instead of distributed control.) The mean value theorem yields
	\begin{align*}
	\frac{1}{t_k}(\overline{B}(S_{\Omega}(D\boldsymbol{v_k})-S_{\Omega}(D\boldsymbol{\overline{v}})),D\boldsymbol{\phi})
	&=
	\int_{\Omega}\overline{B}\int_0^1 S_{\Omega}'\left(D\boldsymbol{\overline{v}}+\theta D(\boldsymbol{v_k}-\boldsymbol{\overline{v}})\right)
	:D\boldsymbol{z_k}:D\boldsymbol{\phi}\, d \theta \, d x
	\\
	&=(\overline{B}S_{\Omega}'(\sigma_k^{\boldsymbol{\phi}})D\boldsymbol{z_k},D\boldsymbol{\phi})
	\end{align*}
	with $\sigma_k^{\boldsymbol{\phi}}:\Omega \to \mathbb{R}$, $\sigma_k^{\boldsymbol{\phi}}(x)=D\boldsymbol{\overline{v}}(x)+\theta_k^{\boldsymbol{\phi}}(x)(D\boldsymbol{v_k}(x)-D\boldsymbol{\overline{v}}(x))$, $\theta_k^{\boldsymbol{\phi}}:\Omega\to [0,1]$. Analogously, we define $\tilde{\sigma}_k^{\boldsymbol{\phi}}$ for the boundary. Then, we can write equation $(\ref{ConvergenceEquation})$ for $\boldsymbol{\phi}:=\boldsymbol{z_k}$ as
	\begin{align}\label{Equality}
	\begin{split}
	&\quad \mu_0(\nabla \boldsymbol{z_k},\nabla \boldsymbol{z_k})+(\overline{B}S_{\Omega}'(\sigma_k^{\boldsymbol{\phi}})D\boldsymbol{z_k},D\boldsymbol{z_k})
	+(\overline{\tau}S_{\Gamma}'(\tilde{\sigma}_k^{\boldsymbol{\phi}})\boldsymbol{z_k},\boldsymbol{z_k})_{\Gamma_b}
	\\
	&=-(BS_{\Omega}(D\boldsymbol{v_k}),\nabla \boldsymbol{z_k})-(\tau S_{\Gamma}(\boldsymbol{v_k}),\boldsymbol{z_k})_{\Gamma_b}.
	\end{split}
	\end{align}
	We want to interpret the left-hand side of equation $(\ref{Equality})$ as a matrix-vector multiplication. We define the indicator function 
	\begin{align*}
	1_A(x)=
	\begin{cases}
	1,\quad x\in A, \\
	0,\quad x\notin A
	\end{cases}
	\end{align*}
	to set
	\begin{align*}
	M^k(x)&:=\left(
	\begin{array}{ccc}
	\mu_0 I & 0 & 0\\
	0 & \overline{B}(x)S_{\Omega}'(\sigma_k^{\boldsymbol{\phi}}(x)) & 0\\
	0 & 0 & 1_{\Gamma_b}(x)\overline{\tau}(x)S_{\Gamma}'(\tilde{\sigma}_k^{\boldsymbol{\phi}}(x)) \\
	\end{array}
	\right),\quad 
	Z_k:=\left(
	\begin{array}{c}
	\nabla \boldsymbol{z_k} \\
	D \boldsymbol{z_k} \\
	\boldsymbol{z_k} \\
	\end{array}
	\right),
	\\
	M(x)&:=\left(
	\begin{array}{ccc}
	\mu_0 I & 0 & 0\\
	0 & \overline{B}(x)S_{\Omega}'(D\boldsymbol{\overline{v}}(x)) & 0\\
	0 & 0 & 1_{\Gamma_b}(x)\overline{\tau}(x)S_{\Gamma}'(\boldsymbol{\overline{v}}(x)) \\
	\end{array}
	\right)
	\end{align*}
	with the identity matrix $I$. We use the standard scalar product on $L:=L^2(\Omega)^{N \times N}\times L^2(\Omega)^{N \times N}\times L^2(\Gamma_b)^N$ for $Z\in L$:
	\begin{align*}
	(Z,Z)_L:=(Z_1,Z_1)_{\Omega}+(Z_2,Z_2)_{\Omega}+(Z_3,Z_3)_{\Gamma_b}
	\end{align*}
	and the induced norm $\|\cdot \|_L$. This yields for the left-hand side of equation $(\ref{Equality})$
	\begin{align*}
	|\mu_0(\nabla \boldsymbol{z_k},\nabla \boldsymbol{z_k})+(\overline{B}S_{\Omega}'(\sigma_k^{\boldsymbol{\phi}})D\boldsymbol{z_k},D\boldsymbol{z_k})+(\overline{\tau}S_{\Gamma}'(\tilde{\sigma}_k^{\boldsymbol{\phi}}) \boldsymbol{z_k},\boldsymbol{z_k})_{\Gamma_b}|
	=
	|(Z_k M^k,Z_k)_{L}|.
	\end{align*}
	Additionally, we use an isomorphism between matrices and vectors by $\nabla \boldsymbol{z_k},D \boldsymbol{z_k}\in H^1(\Omega)^{N \times N}\cong H^1(\Omega)^{N^2}$ to interpret $Z_k$ as a vector. Each block of $M^k$ and $M$ is symmetric and positive-semi definite. Hence, we can find a Cholesky decomposition for each block:
	\begin{align*}
	M^k(x)&=
	\left(
	\begin{array}{ccc}
	M_1(x) & 0 & 0 \\
	0 & M^k_2(x) & 0 \\
	0 & 0 & M^k_3(x) \\
	\end{array}
	\right)
	\\
	&=
	\left(\begin{array}{ccc}
	L_{1}^T(x)L_{1}(x) & 0 & 0 \\
	0 & (L^k_2)^T(x)L^k_2(x) & 0 \\
	0 & 0 & (L^k_3)^T(x)L^k_3(x) \\
	\end{array}\right),
	\\
	M(x)&=		\left(\begin{array}{ccc}
	L_{1}^T(x)L_{1}(x) & 0 & 0 \\
	0 & (L_2)^T(x)L_2(x) & 0 \\
	0 & 0 & (L_3)^T(x)L_3(x) \\
	\end{array}\right)
	\end{align*}
	and set
	\begin{align*}
	L^k(x):=\left(
	\begin{array}{ccc}
	L_1(x) & 0 & 0 \\
	0 & L^k_2(x) & 0 \\
	0 & 0 & L^k_3(x) \\
	\end{array}
	\right),\quad \text{and }
	L(x):=\left(
	\begin{array}{ccc}
	L_1(x) & 0 & 0 \\
	0 & L_2(x) & 0 \\
	0 & 0 & L_3(x) \\
	\end{array}
	\right).
	\end{align*}
	Note that the matrix $L_1$ is independent of $k$. This will make the arguments easier later.
	
	First, we verify $L_2^k\to L_2$ in $L^4(\Omega)^{N \times N}$ and $L_3^k\to L_3$ in $L^4(\Gamma_b)^{N \times N}$ for a subsequence. For that purpose, we show $M^k\to M$ almost everywhere and the boundedness of $M^k$ independent of $k$ to conclude $M_2^k\to M_2$ in $L^2(\Omega)^{N \times N}$ and $M_3^k\to M_3$ in $\times L^2(\Gamma_b)^{N \times N}$. We know
	\begin{align*}
	M_2^k(x)=\overline{B}S_{\Omega}'(\sigma_k^{\boldsymbol{\phi}})=\overline{B}S_{\Omega}'\Big(D\boldsymbol{\overline{v}}(x)+\theta_k^{\boldsymbol{\phi}}(x)(D\boldsymbol{v_k}(x)-D\boldsymbol{\overline{v}}(x))\Big).
	\end{align*}
	The function $\theta_k^{\boldsymbol{\phi}}$ is bounded by one. Lemma $\ref{ControlToStateContinuous}$ implies $\boldsymbol{v_k}\to \boldsymbol{\overline{v}}$ in $V_2$. Thus, there exists an almost everywhere convergent subsequence $\boldsymbol{v_k}\to \boldsymbol{\overline{v}}$. Hence, the continuity of $S_{\Omega}'$ implies
	\begin{align*}
	M_2^k(x)\to M_2(x) \quad \text{almost everywhere.}
	\end{align*}
	The arguments for $M_3^k$ are identical.
	
	Corollary $\ref{SolutionLinearProblem}$ states the existence of a subsequence with $\boldsymbol{z_k}\rightharpoonup \boldsymbol{\overline{z}}$ in $V_2$. Thus, there exists a subsequence with $\boldsymbol{Z_k}\rightharpoonup \boldsymbol{\overline{Z}}$ in $L^2(\Omega)^{N \times N}\times L^2(\Omega)^{N \times N}\times L^2(\Gamma_b)^N$. We conclude for the weak limit $\boldsymbol{\overline{Z}}$ of $\boldsymbol{Z_k}$ with the definition of $L^k$ and equation $(\ref{Equality})$:
	\begin{align*}
	\|L\boldsymbol{\overline{Z}}\|_{L}^2
	&\leq
	\liminf_k \|L^k\boldsymbol{Z_k}\|_{L}^2
	\\
	&=\liminf_k |\mu_0 (\nabla \boldsymbol{z_k},\nabla \boldsymbol{z_k})+(\overline{B}S_{\Omega}'(\sigma_k^{\boldsymbol{\phi}})D \boldsymbol{z_k},D \boldsymbol{z_k})+(\overline{\tau}S_{\Gamma}'(\tilde{\sigma}_k^{\boldsymbol{\phi}})\boldsymbol{z_k},\boldsymbol{z_k})_{\Gamma_b})|
	\\
	&=\liminf_k |(BS_{\Omega}(D\boldsymbol{v_k}),D\boldsymbol{z_k})+(\tau S_{\Gamma}(\boldsymbol{v_k}),\boldsymbol{z_k})_{\Gamma_b}|.
	\end{align*}
	We rewrite the first summand:
	\begin{align*}
	(BS_{\Omega}(D\boldsymbol{v_k}),\nabla \boldsymbol{z_k})=(B(S_{\Omega}(D\boldsymbol{v_k})-S_{\Omega}(D\boldsymbol{\overline{v}})),\nabla \boldsymbol{z_k})+(BS_{\Omega}(D\boldsymbol{\overline{v}}),\nabla \boldsymbol{z_k}).
	\end{align*}
	With $BS_{\Omega}(D\boldsymbol{\overline{v}})\in L^2(\Omega)^{N \times N}$ and $D\boldsymbol{z_k}\rightharpoonup D\boldsymbol{\overline{z}}\in V_2$ follows 
	\begin{align*}
	(BS_{\Omega}(D\boldsymbol{\overline{v}}),\nabla \boldsymbol{z_k})\to (BS_{\Omega}(D\boldsymbol{\overline{v}}),\nabla \boldsymbol{\overline{z}}).
	\end{align*}
	For the first summand, we have with Lemma $\ref{UpperBound}$, and $C,\tilde{C}\in \mathbb{R}$
	\begin{align*}
	&\quad|(B(S_{\Omega}(D\boldsymbol{v_k})-S_{\Omega}(D\boldsymbol{\overline{v}})),\nabla \boldsymbol{z_k})|
	\\
	&\leq C\|B\|_{L^{\infty}(\Omega)}\int_{\Omega}\left(\delta+|D\boldsymbol{v_k}|+|D\boldsymbol{\overline{v}}|\right)^{p-2}|D\boldsymbol{v_k}-D\boldsymbol{\overline{v}}|\, |\nabla \boldsymbol{z_k}|\, d x
	\\
	&\leq
	\tilde{C} \|B\|_{L^{\infty}(\Omega)}\delta^{p-2}\|\boldsymbol{v_k}-\boldsymbol{\overline{v}}\|_{V_2} \|\boldsymbol{z_k}\|_{V_2}\to 0\text{ for }k\to \infty
	\end{align*}
	due to $\boldsymbol{v_k}\to \boldsymbol{\overline{v}}\in V_2$ and the boundedness of $(\boldsymbol{z_k})_k$. 
	Analogously, we have
	\begin{align*}
	(\tau S_{\Gamma}(\boldsymbol{v_k}),\boldsymbol{z_k})_{\Gamma_b}\to (\tau S_{\Gamma}(\boldsymbol{\overline{v}}),\boldsymbol{\overline{z}})_{\Gamma_b}.
	\end{align*}
	Hence, it follows
	\begin{align*}
	\|L\boldsymbol{\overline{Z}}\|_{L}^2\leq \liminf_k\|L^k\boldsymbol{Z_k}\|_{L}^2&=|B(S_{\Omega}(D\boldsymbol{\overline{v}}),D\boldsymbol{\overline{z}})+(\tau S_{\Gamma}(\boldsymbol{\overline{v}}),\boldsymbol{\overline{z}})_{\Gamma_b}|
	\\
	&=|\mu_0(\nabla \boldsymbol{\overline{z}},\nabla \boldsymbol{\overline{z}})+(\overline{B}S_{\Omega}'(\overline{\sigma}^{\boldsymbol{\phi}})D\boldsymbol{\overline{z}},D\boldsymbol{\overline{z}})+(\overline{\tau}S_{\Gamma}'(\overline{\tilde{\sigma}}^{\boldsymbol{\phi}}) \boldsymbol{\overline{z}},\boldsymbol{\overline{z}})_{\Gamma_b}|
	\\
	&=\|L\boldsymbol{\overline{Z}}\|_{L}^2.
	\end{align*}
	The weak convergence and the norm convergence imply $L^k\boldsymbol{Z_k}\to L\boldsymbol{\overline{Z}}$ in $L^2(\Omega)^{N \times N}\times L^2(\Omega)^{N \times N}\times L^2(\Gamma_b)^N$, which implies $L_1 \nabla \boldsymbol{z_k}\to L_1 \nabla \boldsymbol{\overline{z}}$ in $L^2(\Omega)^{N \times N}$. Hence, we find a subsequence, denoted by the index $k$, with 
	\begin{align*}
	L_1\nabla \boldsymbol{z_k}\to L_1\nabla \boldsymbol{\overline{z}}\quad \text{almost everywhere.}
	\end{align*}
	Because $L_1(x)$ is positive definite for $x\in \Omega$, it is invertible. Thereby, we have
	\begin{align*}
	\nabla \boldsymbol{z_k}=L_1^{-1}L_1\nabla \boldsymbol{z_k}\to L_1^{-1}L_1\nabla \boldsymbol{\overline{z}}=\nabla \boldsymbol{\overline{z}}\quad \text{almost everywhere.}
	\end{align*}
	Corollary $\ref{SolutionLinearProblem}$ yields boundedness of a subsequence of $(\boldsymbol{z_k})_k$ in $V_2$, due to the weak convergence. Thus, we can apply dominated convergence and conclude $\nabla \boldsymbol{z_k}\to \nabla \boldsymbol{\overline{z}}$ in $L^2(\Omega)^{N \times N}$.  Due to $\boldsymbol{z_k}\rightharpoonup \boldsymbol{\overline{z}}$ in $V_2$ follows $\boldsymbol{z_k}\to \boldsymbol{\overline{z}}$ in $L^2(\Omega)^N$. Thus, we conclude $\boldsymbol{z_k}\to \boldsymbol{\overline{z}}$ in $V_2$.
\end{proof}
\begin{theorem}\label{GateauxControlToState}
	Let $\mu_0,\delta \in (0,\infty)$. The coefficient-to-state operator $\mathcal{S}$ is Gâteaux differentiable.
\end{theorem}
\begin{proof}
	Let $(\overline{B},\overline{\tau}),(B,\tau)\in W$, $(t_k)_k \in \mathbb{R}^{\mathbb{N}}$ with $t_k\to 0$, $\boldsymbol{\overline{v}}:=\mathcal{S}(\overline{B},\overline{\tau})$, $\boldsymbol{v_k}:=\mathcal{S}(\overline{B}+t_kB,\overline{\tau}+t_k\tau)$, and $\boldsymbol{z_k}:=(\boldsymbol{v_k}-\boldsymbol{\overline{v}})/t_k$.
	
	Lemma $\ref{ControlToStateContinuous}$ yields the boundedness of $(\boldsymbol{z_k})_k\in V_2$. Hence, there exists a weakly convergent subsequence $\boldsymbol{z_{k_j}}\rightharpoonup \boldsymbol{\overline{z}}$. Corollary \ref{SolutionLinearProblem} states that $\boldsymbol{\overline{z}}$ is the solution of
	\begin{align*}
	\langle D_{\boldsymbol{v}}G(\boldsymbol{\overline{v}}(B,\tau))\boldsymbol{\overline{z}},\boldsymbol{\phi}\rangle_{V_2^*,V_2}=-\langle D_{B,\tau}G(\boldsymbol{\overline{v}}(B,\tau),\boldsymbol{\phi}\rangle_{V_2^*,V_2}.
	\end{align*}
	We verified in \cite[Lemma $4.1$]{Schmidt2023a} that there exists only one solution to this problem. Thus, all weakly convergent subsequences converge to $\boldsymbol{\overline{z}}$. It follows weak convergence of the original sequence $\boldsymbol{z_k}\rightharpoonup \boldsymbol{\overline{z}}$ in $V_2$. We proved in Lemma $\ref{StrongConvergence}$ that $\boldsymbol{z_k}\to \boldsymbol{\overline{z}}$ in $V_2$. Therefore, the directional derivative exists. We see from equation $(\ref{Linearity})$ the linearity of $\boldsymbol{v}'(\overline{B},\overline{\tau})(B,\tau)$ because the right-hand side changes linearly and we have a linear problem. Due to the continuity of equation $(\ref{Linearity})$ with respect to $B$ and $\tau$, we have boundedness of
	\begin{align*}
	\sup_{(B,\tau)\in W\text{, }\|(B,\tau)\|_W=1}|\boldsymbol{v'}(\overline{B},\overline{\tau})(B,\tau)|.
	\end{align*}
	Thus, the coefficient-to-state operator $\mathcal{S}$ is Gâteaux differentiable.
\end{proof}
We conclude:
\begin{corollary}
	Let $\mu_0,\delta\in (0,\infty)$. The cost function $f$, see relation $(\ref{OptimalControlProblem})$, defined on $W$ is Gâteaux differentiable.
\end{corollary}
\begin{proof}
	The inner function, the coefficient-to-state operator $\mathcal{S}$, is Gâteaux differentiable, see Theorem $\ref{GateauxControlToState}$. The squared norm is Fréchet differentiable. The trace operator $tr$ and the projection operator $P$ are bounded linear operators. Thus, the directional derivatives are the mappings $(\boldsymbol{v},\boldsymbol{w})\mapsto tr'(\boldsymbol{v})\boldsymbol{w}=tr(\boldsymbol{w})$ and $(\boldsymbol{v},\boldsymbol{w})\mapsto P'(\boldsymbol{v})\boldsymbol{w}=P(\boldsymbol{w})$. These mappings are Lipschitz continuous in both components. Thus, the trace operator and the projection operator are Fréchet differentiable, too. Now, the chain rule as in \cite[Proposition $4.1(b)$]{Zeidler1986} yields the Gâteaux differentiability of the cost function.
\end{proof}
\section{Adjoint equation}\label{DualEquation}
In this section, we formulate the necessary optimality condition for the cost function and prove existence and uniqueness of the solution for the adjoint equation. In general, the necessary optimality condition is:
\begin{lemma}
	Let $W\subseteq H^1(\Omega)\times H^1(\Gamma_b)$ be nonempty and convex, $f:W\to \mathbb{R}$ bounded from below and differentiable in all directions $(B,\tau)\in W$. Let $(\overline{B},\overline{\tau})$ be the local minimum of $f$. Then follows
	\begin{align*}
	f'(\overline{B},\overline{\tau})((B,\tau)-(\overline{B},\overline{\tau}))\geq 0\quad \text{for all }(B,\tau)\in W.
	\end{align*}
\end{lemma}
\begin{proof}
	See \cite[Theorem $3.2$]{DelosReyes2015}.
\end{proof}
As the cost function $f$ is Gâteaux differentiable, the necessary optimality condition for the local minimum $(\overline{B},\overline{\tau})\in W$ is
\begin{align*}
f'(\overline{B},\overline{\tau})((B,\tau)-(\overline{B},\overline{\tau}))\geq 0\quad \text{for all }(B,\tau)\in W.
\end{align*}
However, calculating $f'$ directly is computationally expensive. Instead, one uses the dual equation, \cite[Section $1.6.2$]{Hinze2009}. Up to now, we used the reduced cost function. For using the dual equation, we introduce the cost function depending on the solution and the control $\hat{f}:V_2\times W\to \mathbb{R}$,
\begin{align*}
\hat{f}(\boldsymbol{v},(B,\tau))=\frac{1}{2}\|P(tr\boldsymbol{v})-P(\boldsymbol{v_s})\|_{L^2(\Gamma_a)^N}^2+\frac{\epsilon_1}{2}\|\nabla B\|_{L^2(\Omega)^N}^2+\frac{\epsilon_2}{2}\|\nabla \tau\|_{L^2(\Gamma_b)^N}^2
\end{align*}
and the dual equation: Find $\boldsymbol{\lambda}\in V_2$ with
\begin{align}\label{DualEqu}
\mathcal{B}_{\boldsymbol{v}}(\boldsymbol{\lambda},\boldsymbol{\phi})=- \langle D_{\boldsymbol{v}}\hat{f}(\boldsymbol{v},(B,\tau)),\boldsymbol{\phi}\rangle_{V_2^*,V_2}\quad \text{for all }\boldsymbol{\phi}\in V_2
\end{align}
with $\mathcal{B}_{\boldsymbol{v}}:V_2\times V_2\to \mathbb{R}$,
\begin{align*}
\quad \,\mathcal{B}_{\boldsymbol{v}}(\boldsymbol{\lambda},\boldsymbol{\phi})
&=\int_{\Omega}BS_{\Omega}'(D\boldsymbol{v})D\boldsymbol{\lambda}:\nabla \boldsymbol{\phi}\, d x+\mu_0 \int_{\Omega}\nabla \boldsymbol{\lambda}:\nabla \boldsymbol{\phi}\, d x
+\int_{\Gamma_b}\tau S'_{\Gamma}(\boldsymbol{v})\boldsymbol{\lambda}\cdot \boldsymbol{\phi}\, d s.
\end{align*}
By assuming additionally continuous Fréchet differentiability of $f$ and $A$, and continuous invertability of $D_{\boldsymbol{v}}A$, we obtain for all $(B,\tau)\in W$ with the dual operator $A^*$ of the operator $A$
\begin{align*}
f'(\overline{B},\overline{\tau})
=
D_{B,\tau}A\boldsymbol{v}(\overline{B},\overline{\tau})^*\lambda(B,\tau)+D_{\boldsymbol{v}}\hat{f}(\boldsymbol{v}(B,\tau),(B,\tau))
\end{align*}
with $\lambda$ being the solution of the dual equations, see equation $(\ref{DualEqu})$, for details see \cite[Section $1.6$]{Hinze2009}.
\begin{lemma}
	Let $\mu_0,\delta \in (0,\infty)$. There exists one unique solution $\boldsymbol{\lambda} \in V_2$ for the variational formulation in equation $(\ref{DualEqu})$.
\end{lemma}
\begin{proof}
	We verify continuity and coercivity to apply Lax-Milgram's Theorem, which yields the existence and uniqueness of a solution. Due to the Frécht differentiability of $\hat{f}$, we conclude $D_{\boldsymbol{v}}\hat{f}(\boldsymbol{v}(B,\tau),\boldsymbol{\lambda})\in V_2^*$. Additionally, we have
	\begin{align*}
	|\mathcal{B}_{\boldsymbol{v}}(\boldsymbol{\lambda},\boldsymbol{\phi})|
	&\leq
	\int_{\Omega}\left|B\left(|D\boldsymbol{v}|^2+\delta^2\right)^{(p-2)/2}\left((p-2)\frac{(D\boldsymbol{v}:D\boldsymbol{\lambda})\, (D\boldsymbol{v}:\nabla \boldsymbol{\phi})}{|D\boldsymbol{v}|^2+\delta^2}+D\boldsymbol{\lambda}:\nabla \boldsymbol{\phi}\right)\right|\, d x
	\\
	&\quad + \int_{\Gamma_b}\left| \tau \left(|\boldsymbol{v}|^2+\delta^2\right)^{(s-2)/2}\left(\boldsymbol{\lambda}\cdot \boldsymbol{\phi}+(s-2)\frac{(\boldsymbol{v}\cdot \boldsymbol{\lambda})\, (\boldsymbol{v}\cdot \boldsymbol{\phi})}{|\boldsymbol{v}|^2+\delta^2}\right)\right|\, d s
	+\mu_0 \|\boldsymbol{\lambda}\|_{V_2} \|\boldsymbol{\phi}\|_{V_2}
	\\
	&\leq \|B\|_{L^{\infty}(\Omega)}\delta^{p-2}\int_{\Omega}(2-p)\frac{|D\boldsymbol{v}|^2\, |D\boldsymbol{\lambda}|\, |\nabla \boldsymbol{\phi}|}{|D\boldsymbol{v}|^2+\delta^2}+|D\boldsymbol{\lambda}|\, |\nabla \boldsymbol{\phi}|\, d x
	\\
	&\quad + \|\tau\|_{L^{\infty}(\Gamma_b)}\delta^{s-2}\int_{\Gamma_b}|\boldsymbol{\lambda} |\, |\boldsymbol{\phi}| + (2-s)\frac{|\boldsymbol{v}|^2\, |\boldsymbol{\lambda}|\, |\boldsymbol{\phi}|}{|\boldsymbol{v}|^2+\delta^2}\, ds+\mu_0\|\boldsymbol{\lambda}\|_{V_2}\|\boldsymbol{\phi}\|_{V_2}
	\\
	&\leq \|B\|_{L^{\infty}(\Omega)}\delta^{p-2}\int_{\Omega}(3-p)|D\boldsymbol{\lambda}|\, |\nabla \boldsymbol{\phi}|\, d x
	\\
	&\quad+ \|\tau \|_{L^{\infty}(\Gamma_b)}\delta^{s-2}\int_{\Gamma_b}(3-s)|\boldsymbol{\lambda}|\, |\boldsymbol{\phi}|\, d s
	+ \mu_0 \|\boldsymbol{\lambda}\|_{V_2}\|\boldsymbol{\phi}\|_{V_2}.
	\end{align*}
	The trace operator and Korn's inequality imply the boundedness.
	
	Next, we verify coercivity. With Lemma $\ref{LowBound}$ it follows
	\begin{align*}
	\mathcal{B}_{\boldsymbol{v}}(\boldsymbol{\lambda},\boldsymbol{\lambda})
	\geq \mu_0 \int_{\Omega}|\nabla \boldsymbol{\lambda}|^2\, d x=
	\mu_0 \|\boldsymbol{\lambda} \|_{V_2}^2.
	\end{align*}
\end{proof}
\section{Conclusion}\label{Conclusion}
We showed existence of an optimal solution for the parameter identification problem of the $p$-Stokes equations with an additional diffusion term. The proof has similarities to the proof for the distributed optimal control problem, see \cite{Arada2012}. Similar to \cite{Arada2012}, we can also verify Gâteaux differentiability of the coefficient-to-state operator. We can also formulate the adjoint problem and prove existence and uniqueness of a solution. 

It is an open question if an optimal solution for the parameter identification problem exists for $\mu_0=0=\delta$, and if we have convergence of the optimal solution $\boldsymbol{v}(\mu_0,\delta)$ to $\boldsymbol{v}(0,0)$.
\section*{Acknowledgments}
I thank my supervisor Prof. Thomas Slawig from Kiel University and Prof. Angelika Humbert from the Alfred-Wegener-Institut in Bremerhaven for helpful discussions. Additionally, I thank the reviewers for helpful suggestions to improve the quality of the manuscript.
	\bibliographystyle{fbs}

\end{document}